\documentclass[11pt]{amsart}

\usepackage{amsmath,amscd,amssymb,latexsym}

\newcommand{\sm}{\left(\smallmatrix}
\newcommand{\esm}{\endsmallmatrix\right)}
\newcommand{\mat}{\begin{pmatrix}}
\newcommand{\emat}{\end{pmatrix}}
\renewcommand{\c}{\mathfrak{c}}
\renewcommand{\t}{\tau}
\renewcommand{\l}{\lambda}
\newcommand{\G}{\Gamma}
\renewcommand{\i}{\infty}
\newtheorem{thm}{Theorem}
\newtheorem{lem}[thm]{Lemma}
\newtheorem{cor}[thm]{Corollary}
\newtheorem{prop}[thm]{Proposition}
\numberwithin{equation}{section}
\numberwithin{thm}{section}
\begin{document}

\title[Arithmetic Properties of Traces of singular moduli on congruence subgroups]{Arithmetic Properties of Traces of singular moduli\\ on congruence subgroups}

\author{Soon-Yi Kang and Chang Heon Kim}
\address{Korea Advanced Institute of Science and
Technology, Department of Mathematical Sciences, Daejeon, 305-701
Korea} \email{s2kang@kaist.ac.kr}
\address{Department of Mathematics, Hanyang University,
 Seoul, 133-791 Korea}
\email{chhkim@hanyang.ac.kr}

\begin{abstract}
After Zagier proved that the traces of singular moduli $j(z)$ are Fourier coefficients of a weakly holomorphic modular form, various properties of the traces of the singular values of modular functions mostly on the full modular group $PSL_2(\mathbb{Z})$ have been investigated such as their exact formulas, limiting distribution, duality, and congruences. The purpose of this paper is to generalize these arithmetic properties of traces of singular values of a weakly holomorphic modular function on the full modular group to those on a congruence subgroup $\Gamma_0(N)$.
\end{abstract}
\maketitle

\section{Introduction}

\textit{Singular moduli} are special values of the classical modular invariant $j(z)$ at imaginary quadratic arguments in the upper half plane $\mathbb{H}$.  These important algebraic numbers have a long history in number theory.  Recently, the work of Borcherds \cite{Borch} and Zagier \cite{Zagier} have inspired many works on connecting the traces of singular moduli to the Fourier coefficients of weakly holomorphic modular forms or Harmonic Maass forms of half-integral weight.  (See \cite{ono} for the list of references.)

To describe the context of the work of Borcherds and Zagier, we first define the modular trace of a weakly holomorphic modular function.  In this paper, $D$ is always a positive integer congruent to 0 or 3 modulo 4 unless otherwise specified.  For each non-square $D$,  we let $\mathcal{Q}_{D,N}$ denote the set of positive definite integral binary quadratic forms $$Q(x,y)=[Na,b,c]=Nax^2+bxy+cy^2$$ with discriminant
 $-D=b^2-4Nac$.  The congruence subgroup $\Gamma_0(N)\subseteq \Gamma(1):=PSL_2(\mathbb{Z})$ acts on the set $\mathcal{Q}_{D,N}$.  For a fixed solution $\beta \pmod{2N}$ of $\beta^2 \equiv -D \pmod{4N}$, a smaller set $\mathcal{Q}_{D,N,\beta}
 =\{ [Na,b,c]\in \mathcal{Q}_{D,N} \mid
 b \equiv \beta \pmod{2N} \}$ is still invariant under the action of  $\Gamma_0(N)$.  A significance of these forms is that there is a canonical bijection between $\mathcal Q_{D,N,\beta}/\Gamma_0(N)$ and $\mathcal
Q_{D,1}/\Gamma(1)$ when the discriminant $D$  is not divisible by the square of any prime divisor of $N$ \cite{GKZ}. For each quadratic form $Q$ in $\mathcal Q_{D,N,\beta}$, the corresponding Heegner point on the modular curve $X_0(N)$ is the unique root of $Q(x,1)$,
$$z_Q=\frac{-b+i\sqrt{D}}{2Na}\in \mathbb{H}.$$
Denoting the stabilizers of $Q$ in $\Gamma_0(N)$ by $\Gamma_0(N)_Q$, we define the trace of a weakly holomorphic modular function $f$ on $\Gamma_0(N)$  by
\begin{equation}\label{tracebeta}\textbf{t}_{f}(D) =
 \underset{Q\in\mathcal{Q}_{D,N,\beta}/\Gamma_0(N)}{\sum}
  \frac{1}{|\Gamma_0(N)_Q|} f (z_Q).\end{equation}
In particular, the class number $H_N(D)$ is given by
\begin{equation}H_N(D)=\textbf{t}_{1}(D)=\underset{Q\in\mathcal{Q}_{D,N,\beta}/\Gamma_0(N)}{\sum}
  \frac{1}{|\Gamma_0(N)_Q|}\end{equation}
which is the Hurwitz-Kronecker class number $H(D)$ when the level $N=1$.

We consider a generalized Hilbert class polynomial $\mathcal{H}_D$ by
\begin{equation}\label{cp}\mathcal{H}_D(X)=\prod_{Q\in Q_{D,1}/\Gamma(1)}(X-j(z_Q))^{1/|\Gamma(1)_Q|}\end{equation}
that reduces to the Hilbert class polynomial  when $-D<0$ is a fundamental discriminant. If $q=e^{2\pi i z}$ and $j_1(z)=j(z)-744=q^{-1}+196884q+21493760q^2+\cdots$ denotes the normalized Hauptmodul for $\G(1)$, then the $q$-expansion of the polynomial is given by \cite[Eq.~(11)]{Zagier}
\begin{equation}\label{cpq}\mathcal{H}_D(j(z))=q^{-H(D)}(1-\mathbf{t}_{j_1}(D)q+O(q^2)).\end{equation}
In \cite[p.~204]{Borch}, \cite[Theorem 3]{Zagier}, Borcherds proved that
 \begin{equation}\label{bp}\mathcal{H}_D(j(z))=q^{-H(D)}\prod_{n=1}^\i(1-q^n)^{A(n^2,D)},\end{equation}
where  $A(d,D)$ are the Fourier coefficients of a weight $1/2$ weakly holomorphic modular form.
More precisely, for $D\geq 0$, the set of the functions
$$f_D(z)=q^{-D}+\mathop{\sum_{d >0}}_{d\equiv 0, 3 \, (\textrm{mod} \, 4)}  A(d,D)q^d$$  form a unique basis of the space of weakly holomorphic modular forms of weight $1/2$.  Comparing equations (\ref{cpq}) and (\ref{bp}), we see that $\mathbf{t}_{j_1}(D)=A(1,D)$ for all $D>0$.
On the other hand, Zagier showed that
$$g_d(z)=q^{-d}-\mathop{\sum_{D\geq 0}}_{D\equiv 0, 3 \, (\textrm{mod} \, 4)}  A(d,D)q^D$$ is a weakly holomorphic modular form of weight $3/2$ and the set $\{g_d|d> 0\}$ form a basis of the space of weakly holomorphic modular forms of weight $3/2$.  Accordingly, $\mathbf{t}_{j_1}(D)$ is the coefficient of $q^D$ in the Fourier expansion of $g_1(z)$.

This type of duality relation involving traces of the values of Niebur-Poincar\'{e} series on $\G(1)$ has been found by Bringmann and Ono \cite[Theorems 1.1 and 1.2]{BO}.
For a non-negative integer $m$ and complex numbers $s$ and $z=x+iy$ with $y>0$, we define the weight zero $m$th Niebur-Poincar\'e series on a congruence subgroup $\G$ by
 \begin{equation}\label{poincare}
 \mathcal F_m (z,s)=\sum_{M\in \G_\i\backslash\G}
  e(-m {\rm Re}Mz) ({\rm Im} Mz)^{1/2} I_{s-1/2} (2\pi m {\rm Im}
  Mz),
 \end{equation}
where $\G_\i\subset\G$ is the subgroup of translations, $e(z)=e^{2\pi i z}$, and $I_{s-1/2}$ is the modified Bessel function of the first kind.  Niebur observed that every modular function on $\G(1)$ that is holomorphic away from the cusp at infinity can be written as a linear sum of these Poincar\'{e} series and Duke \cite{Duke} utilized this property of the Poincar\'{e} series to obtain explicit formulas for the traces of CM values of Hecke type Faber polynomials $j_m(z)=P_m(j(z))$ in terms of Kloosterman sums and the class number $H(D)$.  Extending their ideas to a modular function of prime level $p$, the authors and Choi and Jeon \cite{CJKK2} derived  exact formulas for the traces of singular values on $\G_0(p)$.

Moving on to divisibility properties of the modular trace, Ahlgren and Ono \cite{AO} showed that $\textbf{t}_{j_1}(p^2D)\equiv 0 \pmod p$  for an odd prime that splits in $\mathbb{Q}(\sqrt{-D})$ and Osburn \cite{Os} generalized the congruence relation to that of the traces of CM values of Hauptmodul with prime level.  In \cite[Theorem 1.4]{CJKK1}, the authors with Choi and Jeon established Treneer type divisibility property of the traces of singular moduli with arbitrary level.

Although, Bruinier and Funke \cite{BF} generalized Zagier's result by showing that the modular traces of any weakly holomorphic modular functions are Fourier coefficients in the  holomorphic part of a Harmonic weak Maass form, most results on traces of singular moduli are on level 1 case.  In this paper, we generalize the properties of traces of singular moduli discussed in this introduction such as duality, exact formulas, and congruences to the traces of singular values of any modular function that is completely determined by its principal part at infinity.  In the next section, we will state our results with interesting examples. Then in the three sections that follow we give the proofs of our theorems.  In the last section, we discuss very briefly a property of the modular trace that we don't treat in this paper.

\section{statements of results and Examples}\label{results}

\subsection{Exact formulas for traces of singular moduli on congruence subgroups}\label{exact}
Throughout the rest of the paper, $\G$ always denotes the group $\G_0^*(N)$ generated by $\G_0(N)$ and all Atkin-Lehner involutions.  Then the set $\mathcal{Q}_{D,N}$ is invariant under the action of $\G$.  We can define the trace of a weakly holomorphic modular function with respect to $\G$ as

\begin{equation*}
\textbf{t}_{f}^*(D) =
 \underset{Q\in\mathcal{Q}_{D,N}/\Gamma}{\sum}
  \frac{1}{|\Gamma_Q|} f (z_Q).
 \end{equation*}
We denote $H_N^*(D)$ by the corresponding class number, i.e., $\textbf{t}_{1}^*(D)$.  It is easy to see that if $a$ is the number of prime divisors $p$ of $(\beta,N)$ such that
 $p\nmid N/(\beta,N)$, then
 \begin{equation} \label{trace}
 \textbf{t}_{f}^*(D) =
    \frac {1}{2^{a}} \textbf{t}_{f}(D).
 \end{equation}
For this reason, we may study $\textbf{t}_{f}^*(D)$ for $\textbf{t}_{f}(D)$ avoiding the inconvenience of dealing with $\beta$ in $\mathcal{Q}_{D,N,\beta}$.
Furthermore, using a theta lift, we establish the following relation between class numbers that will be needed later.
 \begin{thm} \label{cn}
For every positive integer $D$ congruent to $0$ or $3$ modulo 4,
\begin{equation}\label{hh} H_{4}^*(4D)=3H(D). \end{equation}
\end{thm}

Recall the Niebur-Poincar\'{e} series $\mathcal F_m (z,s)$ on $\G$ in (\ref{poincare}). When $m=0$, the series is the non-holomorphic Eisenstein series for $\G$ of weight zero.  The Poincar\'e series $\mathcal F_m (z,s)$ converges absolutely for ${\rm Re} \, s > 1$ and can be analytically continued to the entire $s$ plane and it has no poles at $\mathrm{Re}(s)=1$ \cite[Theorem 5]{Niebur}.  Moreover, it is an eigenfunction for the hyperbolic Laplacian
   $\triangle=-y^2(\partial^2_x+\partial^2_y)$ with eigenvalue $s(1-s)$ so that $\mathcal F_m (z,1)$ is annihilated by the Laplacian $\triangle$ and thus almost holomorphic on $\mathbb{H}$.  Niebur showed that any modular function on $\G(1)$ that is holomorphic away from the cusp at infinity is a linear combination of the  $\mathcal F_m (z,1)$ on $\G(1)$ \cite[Theorem 6]{Niebur}.  Using the similar argument in \cite[p.~4]{CJKK2}, we can generalize the result to $\mathcal F_m (z,1)$ on $\G$.
 \begin{prop}\label{triv} Suppose $f$ is a modular function on $\G$ whose poles are supported only at infinity and the principal part is given by $\sum_{m=1}^\ell a_m e(-mz)$. If we define $$\mathcal F_m^*(z,s)=(2\pi  \sqrt{m})\mathcal F_m
 (z,s)+c_m,$$ where $-c_m$ is the constant term in $(2\pi
\sqrt{m})\mathcal F_m (z,1)$, then
 \begin{equation*}
 f(z)=\sum_{m=1}^\ell a_m \mathcal F_m^*(z,1).
 \end{equation*}\end{prop}

 For example, a Faber polynomial $j_{m}(z)=m(j-744)|_{T_m}$ on $\Gamma(1)$ satisfies
 \begin{equation}\label{jn}j_{m}(z)=\mathcal F_m^*(z,1)=2\pi \sqrt m \mathcal{F}_{m}(z,1)-24\sigma(m),\end{equation}
 where $|_{T_m}$ is the weight $0$ Hecke operator and $\sigma(m)$ is the divisor function \cite{BJO, Duke}.
If $N$ is a prime $p$ and $p^\alpha|| m$, then for a modular function $f$ satisfying the conditions in Proposition \ref{triv}, we have \cite{CJKK2}
\begin{equation}\label{fp}
 f(z)=\sum_{m=1}^\ell a_m \left(2\pi \sqrt{m} \mathcal F_m(z,1)-24\left( \frac{-p^{\alpha+1}}{p+1} \sigma(m/p^\alpha)
 +\sigma(m)\right)\right).
 \end{equation}

The constant $c_m$ for arbitrary level $N$ can be determined explicitly using properties of Ramanujan sum.

\begin{thm}\label{cons} Let $\mathcal F_m (z,s)$ be the Poincar\'e series defined in (\ref{poincare}).
Then the constant term $-c_m$ in $(2\pi
\sqrt{m})\mathcal F_m (z,1)$ is given by
\begin{equation}\label{cm}
-c_m=24\frac{m}{m_N}\sigma(m_N)\prod_{p\mid N}(1-p^{-2})^{-1}  \sum_{e\|N}\frac{1}{e}\prod_{p|\frac{N}{e}} \delta_e(p)p^{1-\beta_{e,p}}\left((1-p^{\beta_{e,p}-\alpha_p-2})(1+p^{-1})-1\right).\end{equation}

Here $\alpha_p:=ord_p(m)$, $\beta_{e,p}=ord_p(\frac{N}{e})$, $m_N$ denotes the largest exact divisor of $m$ satisfying $(m_N,N)=1$, and
$\delta_e(p)$ is defined by
\begin{equation}\label{delta}\delta_e(p)=\left\{%
\begin{array}{ll}
    1, & \hbox{ if $\beta_{e,p}\leq \alpha_p+1$,} \\
    0, & \hbox{ otherwise}. \\
    \end{array}%
\right.\end{equation}
\end{thm}

The constant $c_m=-24\sigma(m)$ when $N=1$  and $c_m=-24\left( \frac{-p^{\alpha+1}}{p+1} \sigma(m/p^\alpha)
 +\sigma(m)\right)$ when $N=p$ and $p^\alpha|| m$.  (See (\ref{jn}) and (\ref{fp}), respectively.)  If $N=p^t$, a prime power, then
\begin{equation}\label{cmt}c_m=-24\frac{p^\alpha\sigma(m/p^\alpha)}{1-p^{-2}}\times
\left\{%
\begin{array}{ll}
    p^{1-t}((1-p^{t-\alpha-2})(1+p^{-1})-1)+p^{-t}, & \hbox{ if $t\leq \alpha+1$,} \\
    p^{-t},
    & \hbox{ otherwise.} \\
\end{array}%
\right.\end{equation}
\par\noindent
Moreover, if $m$ is relatively prime to $N$, then
 \begin{equation}\label{rp}
 c_m= -24\sigma (m) \frac{\sum_{e||N}\mu(N/e)e}{N^2 \prod_{p|N} (1-p^{-2})},
 \end{equation}
 where $\mu(n)$ is the M\"{o}bius function.  This implies  that
 \begin{equation}\label{cm1}
 c_m=\sigma(m)c_1 \qquad \mathrm{whenever}\qquad (m,N)=1.\end{equation}

By Proposition \ref{triv}, one may write the traces of the singular values of $f$ in terms of the traces of Poincar\'e series $\mathcal F^*_m (z,1)$ that lead to exact formulas for the traces of the values of $f$ at Heegner points.  We first find explicit formulas for the traces of the values of $\mathcal F^*_m (z,1)$ at Heegner points.

\begin{lem}\label{trF} Let $\mathcal F_m^*(z,s)=(2\pi  \sqrt{m})\mathcal F_m
 (z,s)+c_m,$ where $\mathcal F_m (z,s)$ is defined in (\ref{poincare})
 and $-c_m$ is the constant term in $(2\pi
\sqrt{m})\mathcal F_m (z,1)$ which is explicitly given in (\ref{cm}). Then the trace of the values of $\mathcal{F}_m^*$ at Heegner points is
 given by
\begin{equation}\label{trm}\textbf{t}_m(D):=\sum_{Q\in\mathcal Q_{D,N}/\Gamma} \frac{1}{|\Gamma_{Q}|}
 \mathcal F_m^*(z_{Q},1)= c_m H_{N}^*(D)+
 \mathop{\sum_{c>0}}_{c\equiv 0 \, (\textrm{mod} \, 4N)} S_D(m,c)\sinh\left(\frac{4\pi m \sqrt
 D}{c}\right),\end{equation}
where
\begin{equation}\label{salie}S_D(m,c)=\sum_{x^2\equiv -D (\textrm{mod} \, c)} e(2mx/c)\quad \mathrm{for\ any\ positive\ integers\ }m\ \mathrm{and}\ c.\end{equation}
\end{lem}

It follows from Proposition \ref{triv} and Lemma \ref{trF} that

\begin{cor}\label{main}
Suppose that $f$ is a modular function for $\G$
whose poles are supported only at $\i$ and the principal part is given by $\sum_{m=1}^\ell a_m e(-mz)$.
Then \begin{equation}\label{ef}
 \textbf{t}_f^*(D)=\sum_{m=1}^\ell a_m
 \left[
 c_m H_{N}^*(D)+\mathop{\sum_{c>0}}_{c\equiv 0\,(\textrm{mod}\, 4{N})}
 S_D(m,c)\sinh\left(\frac{4\pi m\sqrt{D}}{c}\right)
 \right],
 \end{equation}
 where $-c_m$ is given in (\ref{cm}).\end{cor}
\noindent {\bf Example 1.}
Consider $$j_{45}^*=-1+\left(\frac{\eta(3z)^2 \eta(15z)^2}{\eta(z)\eta(5z)\eta(9z)\eta(45z)}\right),$$
where $\eta(z)$ is the Dedekind eta function defined by $\eta(z)=q^{1/24} \prod_{n=1}^{\infty} (1-q^n)$. Then $j_{45}^*$ is the Hauptmodul for $\Gamma_0^*(45)$ which is of genus $0$ and has a Fourier expansion of the form $q^{-1}+0+O(q)$. Let $D=20$. Since the representatives for $\mathcal Q_{20,45,40}/\Gamma_0(45)$ are given by
$[405,40,1]$ and $[90,-50,7]$, we find from equations (\ref{trace}), (\ref{rp}) and (\ref{ef}) that
$$ j_{45}^*(z_{[405,40,1]})+j_{45}^*(z_{[90,-50,7]})
=-1 + 2\mathop{\sum_{c>0}}_{c\equiv 0\,(\textrm{mod}\, 180)}
 S_{20}(1,c)\sinh\left(\frac{8\pi \sqrt{5}}{c}\right).
$$
The left hand side of the equation is known to be $-3$, and hence the exponential sum on the right hand side has the value $-1$.
\bigskip

\noindent {\bf Example 2.} Let $P_m(j_4^*)$ denote the polynomial of the Hauptmodul $j_4^*(z)$ on $\Gamma_0^*(4)$
with Fourier development of the form $P_m(j_4^*)(z)=q^{-m} + O(q)$.  For an odd integer $m$, we obtain from Corollary \ref{main} and (\ref{rp}) that
\begin{equation}\label{41} \textbf{t}_{P_m(j_4^*)}^*(4D)=-8\sigma (m) H_{4}^*(4D)+\mathop{\sum_{c>0}}_{c\equiv 0\,(\textrm{mod}\, 16)}
 S_{4D}(m,c)\sinh\left(\frac{4\pi m\sqrt{4D}}{c}\right). \end{equation}
But by an analogous proof to \cite[Theorem 1.2]{CJKK1}, we find that
\begin{equation}\label{42} \textbf{t}_{P_m(j_4^*)}^*(4D) = -24 \sigma (m) H(D). \end{equation}
Comparing two equations (\ref{41}) and (\ref{42}) and applying Theorem \ref{cn}, we find that
\begin{equation} \mathop{\sum_{c>0}}_{c\equiv 0\,(\textrm{mod}\, 16)}
 S_{4D}(m,c)\sinh\left(\frac{4\pi m\sqrt{4D}}{c}\right)=0
\end{equation}
 for all positive integer $D\equiv 0,3 \,(\textrm{mod}\, 4)$.

\subsection{Duality and Divisibility}
For more arithmetic properties of traces of the values of Niebur-Poincar\`{e} series and modular functions holomorphic away from the cusp at infinity, we make the Kloosterman sum representation of the trace of the Niebur-Poincar\`{e} series.  In order to define Kloosterman sum, we need the extended Kronecker symbol $\displaystyle{\left(\frac{c}{d}\right)}$ and
$$\varepsilon_d:=\left\{
                     \begin{array}{ll}
                       1, & \hbox{ if $d\equiv 1 \pmod 4$,} \\
                       i, & \hbox{ if $d\equiv 3 \pmod 4$}
                     \end{array}
                   \right.$$
that is defined for odd $d$.
For $c, m, n, \l \in \mathbb{Z}$ with $c \equiv 0 \pmod 4$, the weight $k:=\l+1/2$ Kloosterman sum $K_\lambda(m,n,c)$ is defined by
\begin{equation}\label{kloos}
K_\lambda(m,n,c):=\sum_{v \pmod c^*}\left(\frac{c}{v}\right)\varepsilon_v^{2\lambda+1}e\left(\frac{m\bar{v}+nv}{c}\right),\end{equation}
where the sum runs through the primitive residue classes modulo $c$ and $v\bar{v}\equiv 1 \pmod c$.

Following the method developed in \cite{Hej}, \cite{Bru}, \cite{BJO}, \cite{BO}, and \cite{MP}, we construct a half integral weight Maass-Poincar\'e series
for arbitrary level $4N$ whose holomorphic coefficients are represented by the Kloosterman sums:  For $s\in \mathbb{C}$ and $y\in \mathbb{R}-\{0\}$, we define
$$
\mathcal{M}_s(y):=|y|^{-k/2}M_{\frac{k}{2}\mathrm{sgn}(y),s-\frac{1}{2}}(|y|),$$
where $M_{\nu,\mu}$ is the usual M-Whittaker function.  And for $m\geq 1$ with $(-1)^{\l+1}m\equiv 0,1 \pmod 4$, we define
$$\varphi_{-m,s}(z):=\mathcal{M}_s(-4\pi m y) e(-mx).$$
With these notations, we define the Poincar\'e series for $\mathrm{Re}(s)>1$ by

\begin{equation}\label{poinhalf}
\mathfrak{P}_{\l,n}(-m,s;z):=\sum_{M\in \G_\i\backslash\G_0(4N)}(\varphi_{-m,s}|_k M)(z)\end{equation}
 where  $|_k$ is the usual weight $k$ slash operator.

Now as in \cite{BO} and \cite{MP}, we apply Kohnen's projection operator \cite[p.~250]{Kohnen} $\mathrm{pr}_\l$ to (\ref{poinhalf}) to obtain a new family of weak Maass forms.

\begin{equation}
\mathcal{P}_{\l,N}(-m,z):=\left\{
                  \begin{array}{ll}
                    \frac{3}{2}\mathfrak{P}_{\l,n}(-m,\frac{k}{2};z)|\mathrm{pr}_\l, & \hbox{if $\l\geq 1$;} \\
                    \frac{3}{2(1-k)\G(1-k)}\mathfrak{P}_{\l,n}(-m,1-\frac{k}{2};z)|\mathrm{pr}_\l, & \hbox{if $\l\leq 0$.}
                  \end{array}
                \right.
\end{equation}

 The series $\mathcal{P}_{\l,N}(-m,z)$ is a weakly holomorphic modular form of weight $\l+1/2$ and level $4N$ satisfying Kohnen plus-condition if $\l>1$ and it is a weak Maass form if $\l\leq 1$ that has Fourier expansion
\begin{equation}\mathcal{P}_{\l,N}(-m,z)=q^{-m}+\mathop{\sum_{n\geq 0}}_{(-1)^\l n\equiv 0,1 \pmod 4}b_{\l,N}(-m,n)q^n+\mathcal{P}^{-}_{\l,N}(-m,z),\end{equation}
where $\mathcal{P}^{-}_{\l,N}(-m,z)$ is the non-holomorphic part.

It follows from \cite[Theorem 2.1]{MP} that if $m,n$ are positive integers such that $(-1)^{\l+1} m, (-1)^\l n\equiv 0,1 \pmod 4$ and $N$ is odd, then the Fourier coefficients $b_{\l, N}(-m,n)$  of the weak Maass form $\mathcal{P}_{\l,N}(-m,z)$ is given by
\begin{equation}\label{mkloo}
b_{\l, N}(-m,n):=(-1)^{[\frac{\l+1}{2}]}\pi \sqrt 2(n/m)^{\frac{2\l-1}{4}}(1-(-1)^\l i)\mathop{\sum_{c>0}}_{4N|c}\frac{K_\lambda(-m,n,c)}{c}\delta_o(c/4)I_{\l-1/2}\left( \frac{4\pi\sqrt{nm}}{c} \right),
\end{equation}
where $\delta_o(d)=2$ if $d$ is odd and $1$ otherwise.

Now, the trace of the values of Niebur-Poincar\'e series at Heegner points is  a linear sum of these coefficients and class numbers.
\begin{lem} \label{Kloo}  Let $\textbf{t}_m(D)$ be the trace of the Poincar\'e series of odd level $N$ at Heegner points in (\ref{trm}).
If $(m,N)=1$, then
$$\textbf{t}_m(D)=-\sum_{\nu\mid m}\nu B(-\nu^2,D),$$
where $B(-m,n)$ is given by
\begin{equation}\label{B}
B(-m,n)=-c_1\delta_{\Box}(m)H_N^*(n)+b_{1, N}(-m,n).\end{equation}
Here $-c_1$ is the constant given in (\ref{rp}) with $m=1$, the function $\delta_{\Box}(m)=1$ if $m$ is a square and zero otherwise, and $b_{1, N}(-m,n)$ is given in (\ref{mkloo}) with $\l=1$.
\end{lem}

Generalizing the duality relation $b_{\l,1}(-m,n)=-b_{1-\l,1}(-n,m)$ by Bringmann and Ono \cite[Theorem 1.1]{BO}, we can establish a duality relation for $B(-m,n)$ in  (\ref{B}).  Using the Bruinier-Funke theta lift as in \cite{CJKK1}, one can construct a weight $3/2$ Harmonic weak Maass form $G_N (z)$ whose
 holomorphic part is $\sum_D H_N^* (D) q^D$. Let
 \begin{equation}\label{p1}\mathcal{P}_{1,N}^* (-m,z):=\mathcal{P} _{1,N}(-m,z)+(-c_1)\delta_{\square}(m) G_N (z).\end{equation}
  In addtion, we define
 \begin{equation}\label{p0}\mathcal{P}_{0,N}^* (-m,z):= \mathcal{P}_{0,N}(-m,z)+c_1H_N^*(m)\theta (z)/2,\end{equation}
where $\theta(z)=\sum_{n\in \mathbb{Z}} q^{n^2}$ is the Jacobi theta series.  Then we have the following duality relation between coefficients of $\mathcal{P}_{1,N}^* (-m,z)$ and $\mathcal{P}_{0,N}^* (-m,z)$.
\begin{thm}\label{duality}
Let $\mathcal{P}_{1,N}^* (-m,z)$ and $\mathcal{P}_{0,N}^* (-m,z)$ be weak Maass forms defined in (\ref{p1}) and (\ref{p0}), respectively, with $N$ odd.  Assume they have Fourier coefficients of $q^n$ for $n\geq 0$, $B_{1,N}(-m,n)$ and $B_{0,N}(-m,n)$ , respectively.
If $m$ is a positive integer that is a square modulo ${4N}$, then for every positive integer $n$ with $-n\equiv \square \pmod {4N}$, we have
$$ B_{1,N}(-m,n)=-B_{0,N}(-n,m). $$
\end{thm}
Note that $B(-m,n)$ in (\ref{B}) is the coefficient $ B_{1,N}(-m,n)$ in the holomorphic part of $\mathcal{P}_{1,N}^* (-m,z)$.

Lastly, we examine a divisibility property of the traces of singular moduli that is a generalization of \cite{AO} and \cite{Os}.
\begin{thm}\label{last}  Suppose that $f$ is a modular function for $\Gamma$ with odd level $N$ which is holomorphic away from the cusp at infinity.
If $f$ has the principal part $\sum_{m=1}^\ell a_{m} e(-mz)$ at $\infty$ with $a_{m}\in\mathbb Z$ such that for all non-zero $a_m$, $(m,N)=1$, then  for every prime $p$ for which $(p,N)=1$ and $\left(\frac{-D}{p}\right)=1$,
\begin{equation} \label{cons1}
\textbf{t}_{f}(p^2 D)=-p\sum_{m=1}^\ell a_m\left(\sum_{\nu\mid m}\nu B(-p^2\nu^2,D)+\sum_{\nu\mid m/p}\nu B(-\nu^2,D)\right).
\end{equation}
\end{thm}

If we assume  a modular function $f\in \mathbb Q ((q))$ satisfies the condition in Theorem \ref{last} and both trace of $f$ and the sum $\sum_ma_m(\sum_{\nu\mid m}\nu B(-p^2\nu^2,D)+\sum_{\nu\mid m/p}\nu B(-\nu^2,D))$ are integers, we have the congruence $$\textbf{t}_f(p^2D)\equiv 0 \pmod p.$$

\noindent {\bf Example 3.} Consider
 $$f=\left( \frac{\eta(z)}{\eta(37z)} \right)^2-2
 +37\left( \frac{\eta(37z)}{\eta(z)}\right)^2.$$
 Then $f$ is a modular function for $\Gamma_0^*(37)$ which is of genus 1 and has a
 Fourier expansion of the form
 $ q^{-3}-2q^{-2}-q^{-1}+0+O(q)$.  When $D=11$, the computation from Theorem \ref{cons}, Lemma \ref{trF} and Corollary \ref{main} shows that $\textbf{t}_f(11)= \textbf{t}_3(11)
  -2 \textbf{t}_2(11) - \textbf{t}_1(11) = 5 $
is an integer while $\textbf{t}_1(11)$,  $\textbf{t}_2(11)$,
  $\textbf{t}_3(11)$ are not integers.
Furthermore, $\textbf{t}_f(p^2\cdot 11)\equiv 0 \pmod p$ for a prime $p$ for which $\left(\frac{-11}{p}\right)=1$.  For example, for the first three of such primes $3$, $5$ and $23$, we have $\textbf{t}_f(99)=\textbf{t}_f(3^2\cdot 11)=-6 \equiv 0 \pmod 3$,
 $\textbf{t}_f(275)=\textbf{t}_f(5^2\cdot 11)=-75 \equiv 0 \pmod 5$, and
 $\textbf{t}_f(5819)=\textbf{t}_f(23^2\cdot 11)=246920364 \equiv 0 \pmod{23}$.

\section{Exact Formulas for traces}
In this section, we will give the proofs of theorems listed in Section \ref{exact}.

\begin{proof}[Proof of Proposition \ref{triv}] First, we recall that the Niebur-Poincar\'e series $\mathcal F_m (z,s)$ in (\ref{poincare}) has the following Fourier expansion
 \cite[Theorem 1]{Niebur}; for ${\rm Re} \, s > 1$,
 \begin{equation}\label{fourier}
 \mathcal F_m (z,s)=e(-mx)y^{1/2} I_{s-1/2}(2\pi m y)
  + \sum_{n=-\infty}^{\infty} b_n (y,s;-m) e(nx),
 \end{equation}
 where $b_n (y,s;-m) \to 0$ ($n\neq 0$) exponentially as $y\to
 \infty$. Hence the pole of $\mathcal F_m (z,1)$ at infinity
may occur only in
 $e(-mx)y^{1/2} I_{1/2}(2\pi m y)$, which is equal to
 \begin{equation} \label{pole}
 \frac{1}{\pi y^{1/2} m^{1/2}} \sinh (2\pi m y) y^{1/2} e(-mx)
 =\frac{1}{2\pi  \sqrt{m}} (e(-mz)-e(-m\bar z)).
 \end{equation}
 We normalize $\mathcal F_m (z,1)$ by  multiplying with $2\pi
 \sqrt{m}$ so that the coefficient of $e(-mz)$ is normalized.
It is easy to check that $$\mathcal F_m^*(z,1)=(2\pi  \sqrt{m})\mathcal F_m
 (z,1)+c_m,$$ where $-c_m$ is the constant term in $(2\pi
\sqrt{m})\mathcal F_m (z,1)$, is a $\Gamma$-invariant harmonic function and
 $\mathcal F_m^*(z,1)-e(-mz)$ has a zero at $\infty$.
 Hence Proposition \ref{triv} follows from \cite[Theorem 6]{Niebur}.\end{proof}

As the explicit formulas for $\textbf{t}_m(D)$ given in Lemma \ref{trF} can be derived in a very similar way to Lemma 3 in \cite{CJKK2}, we omit the proof of Lemma \ref{trF}.  Hence in order to complete the proof of exact formulas for the traces of the values of a weakly holomorphic modular function on $\G$, it remains to evaluate the constant term $c_m$ given in Theorem \ref{cons}.

\begin{proof}[Proof of Theorem \ref{cons}] It follows from (\ref{fourier}) and \cite[Theorem 1]{Niebur} that the constant term in $(2\pi  \sqrt{m})\mathcal F_m
(z,1)$ is
\begin{eqnarray}\label{cmf}-c_m&=&\lim_{s\to 1} 2\pi  \sqrt{m} b_0(y,s,-m)\cr
&=& \lim_{s\to 1} 2\pi \sqrt{m} (2\pi^s m^{s-1/2}
 \phi_m(s)/\Gamma(s)) y^{1-s}/(2s-1)
 =4\pi^2 m \lim_{s\to 1} \phi_m(s). \end{eqnarray}
Here $\phi_m(s)=\sum_{c>0} S(m,0;c) c^{-2s}$
 and $S(m,n;c)$ is the general Kloosterman sum
 $\sum_{0\le d < |c|} e((ma+nd)/c)$ for $\sm a & * \\ c&d \esm\in
 \Gamma$.  As $M \in \Gamma_0^*(N)$ if and only if $M$ is of the form
 $\mat \sqrt{e} \mathfrak{a }& \mathfrak{b}/\sqrt{e} \\ N\mathfrak{c}/\sqrt{e} &  \sqrt{e} \mathfrak{d} \emat_{\det =1}$
 for some $e\|N$ with $\mathfrak{a,b,c,d}\in \mathbb Z$ and $(\mathfrak{a}, N\c/e)=1$, we can identify the sum
$S(m,0;c)=\sum_{0\le d < |c|} e(ma/c)$ with the sum of $m$-th powers of primitive $N\mathfrak{c}/e$-th
 roots of unity. Hence if we denote the Ramanujan sum by $u_n(q)$, that is the sum of $n$-th powers of primitive $q$-th
 roots of unity,  we find that
 \begin{eqnarray}\label{phi}\phi_m(s) &=& \sum_{c>0} S(m,0;c) c^{-2s}=\sum_{e\|N} \mathop{\sum_{\c\geq 1}}_{(e,\c)=1}u_m(N\c/e)\frac{N^{-2s}}{e^{-s}} \c^{-2s}\cr
 &=& \sum_{e\|N}e^{-s} \mathop{\sum_{\c\geq 1}}_{(e,\c)=1}u_m(N\c/e)(N\c/e)^{-2s}\\
 &=& \sum_{e\|N}e^{-s} \left(\prod_{p|\frac{N}{e}}\sum_{k\geq ord_p(\frac{N}{e})}u_m(p^k)p^{-2sk}\prod_{p\nmid N}\sum_{k\geq 0}u_m(p^k)p^{-2sk}\right),\nonumber
\end{eqnarray}
where the last equality holds due to the multiplicative property of $u_n(q)$ as a function of $q$ and the fact $(e,\c)=1$.  Recall the following known fact on the Ramanujan sum:
\begin{equation}\label{rsv}u_n(p^k)=\left\{%
\begin{array}{ll}
    0, & \hbox{ if $p^{k-1}\nmid n$,} \\
    -p^{k-1}, & \hbox{ if $p^{k-1} || n$,} \\
    \varphi(p^k), & \hbox{ if $p^k \mid n$,} \\
\end{array}%
\right.\end{equation}
where $\varphi(n)$ is Euler's totient function.  Letting $\alpha=\alpha_p=ord_p(m)$ and $m_N$ be the largest exact divisor of $m$ satisfying $(m_N,N)=1$, we deduce from (\ref{rsv}) that
\begin{eqnarray}\label{prod2}
\prod_{p\nmid N}\sum_{k\geq 0}u_m(p^k)p^{-2sk}&=&\prod_{p\nmid N}(1+\varphi(p)p^{-2s}+\varphi(p^2)p^{-4s}+\cdots+\varphi(p^\alpha)p^{-2\alpha s}-p^\alpha p^{-2(\alpha+1)s})\cr
&=&\prod_{p\nmid N} (1-p^{-2s})(1+p^{1-2s}+\cdots+(p^{1-2s})^\alpha)\\
&=& \frac{\zeta(2s)^{-1}}{\prod_{p\mid N}(1-p^{-2s})}\sigma_{1-2s}(m_N).\nonumber\end{eqnarray}
Also, if $\beta=\beta_{e,p}=ord_p(\frac{N}{e})$ and $\delta_e(p)$ is defined as in (\ref{delta}), then we obtain from (\ref{rsv}) that
\begin{eqnarray}\label{prod1}
\prod_{p|\frac{N}{e}}\sum_{k\geq \beta}u_m(p^k)p^{-2sk}&=&\prod_{p|\frac{N}{e}}\delta_e(p)\left(
\varphi(p^\beta)p^{-2\beta s}+\varphi(p^{\beta+1})p^{-2(\beta+1)s}+\cdots+\varphi(p^\alpha)p^{-2\alpha s}-p^\alpha p^{-2(\alpha+1)s}\right)\cr
&=& \prod_{p|\frac{N}{e}}\delta_e(p)
(p^{1-2s})^{\beta-1}\left[(1+p^{1-2s}+\cdots+(p^{1-2s})^{\alpha-\beta+1})(1-p^{-2s})-1\right]\\
&=& \prod_{p|\frac{N}{e}}\delta_e(p)(p^{1-2s})^{\beta-1}\left(\frac{1-(p^{1-2s})^{\alpha-\beta+2}}{1-p^{1-2s}}(1-p^{-2s})-1\right).\nonumber\end{eqnarray}

Therefore, the theorem follows from  (\ref{cmf}), (\ref{phi}), (\ref{prod2}), and (\ref{prod1}).
\end{proof}

We close this section with the proof of Theorem \ref{cn}.

\begin{proof}[Proof of Theorem \ref{cn}] Following the notations in \cite{BF, CJKK1} we denote
the theta lifting of $f$ by $I(\tau, f)$ and Zagier's Eisenstein series of weight $3/2$ by $F(\tau)$ which is given by
\begin{equation}\label{eis}
F(\t)=\sum_DH(D)q^D+\frac{1}{16\pi\sqrt{v}}\sum_{m\in \mathbb{Z}}\beta(4\pi v m^2)q^{-m^2},\end{equation}
where $\t=u+iv$ with $v>0$ and $\beta(s)=\int_1^\i t^{-3/2}e^{-st}dt$.  Bruinier and Funke \cite{BF} showed that the holomorphic part of $I(\tau, f)$ is the generating function of traces of the values at CM points of $f$.  A little careful application of this result as in \cite[eq. (4.3)]{CJKK1} shows that the holomorphic part of $I(\tau,1)|_{U_4}$ is the generating function of class numbers $4H_4^*(4D)$, whereas the non-holomorphic part of
 $I(\tau,1)|_{U_4}$ is given from \cite[(3.5) and Section 4]{CJKK1} by
 $$
 \frac{3\cdot 1 \cdot \frac 12}{4\pi \sqrt {v/4}} \sum_{m\in \mathbb Z} \beta(4\pi v m^2) q^{-m^2}.
 $$
This implies that the difference between
 $I(\tau,1)|_{U_4}$ and $12 F(\tau)$ lies in the space of holomorphic modular forms of weight 3/2
 on $\Gamma_0(4)$ satisfying Kohnen's plus condition, which should be zero.
 Thus we have the identity $I(\tau,1)|_{U_4}=12 F(\tau)$. Comparing the coefficients in $q^D$, the assertion follows.
\end{proof}

\section{Kloosterman sum representation and duality}
We first prove the Kloosterman sum representation of the trace of the values of Niebur-Poincar\'e series at Heegner points in Lemma \ref{Kloo} and then the duality relation given in Theorem \ref{duality}.

\begin{proof}[Proof of Lemma \ref{Kloo}]
From Kohnen's result on the relation between Kloosterman sum and Sali\'{e} sum given in \cite[Proposition 5]{Kohnen} and \cite[Proposition 2.2]{MP}, we can write
the sum in the far right side in (\ref{trm}) as
\begin{eqnarray}\label{sk}\mathop{\sum_{c>0}}_{c\equiv 0 \, (\textrm{mod} \, 4N)} S_D(m,c)\sinh\left(\frac{4\pi m \sqrt
 D}{c}\right)&=&\mathop{\sum_{c>0}}_{4N|c}\sum_{d|(\frac{c}{4},m)}(1+i)(c/d)^{-1/2}\delta_{o}(c/4d)\cr
&\times &K_1(-\frac{m^2}{d^2},D,c/d)\left(\frac{2m\pi^2\sqrt D}{c}\right)^{1/2}I_{1/2}\left(\frac{4\pi m \sqrt D}{c}\right).\end{eqnarray}
Since $(m,N)=1$ and $c_m=\sigma(m)c_1$ when $(m,N)=1$, we deduce from (\ref{trm}) and (\ref{sk}) that
\begin{eqnarray}\textbf{t}_m(D)=c_1 H_{N}^*(D)\sum_{\nu\mid m}\nu&+&\sum_{\nu\mid m}\mathop{\sum_{c>0}}_{4\nu N|c}(1+i)(c/\nu)^{-1/2}\delta_{o}(c/4\nu)\cr
&\times &K_1(-\frac{m^2}{\nu^2},D,c/\nu)\left(\frac{2m\pi^2\sqrt D/\nu}{c/\nu}\right)^{1/2}I_{1/2}\left(\frac{4\pi m \sqrt D/\nu}{c/\nu}\right).\end{eqnarray}
Replacing $c/\nu$ by $c$ and $m/\nu$ by $\nu$, we complete the proof.
\end{proof}

 Next, we prove the duality between the coefficients of $\mathcal{P}_{1,N}^* (-m,z)$ and $\mathcal{P}_{0,N}^* (-m,z)$ in (\ref{p1}) and (\ref{p0}).
 \begin{proof}[Proof of Theorem \ref{duality}]
According to the Fourier development of $\mathcal{P}_{\lambda,N}(-m,z)$ computed in \cite[Theorem 2.1]{MP}, the coefficient $b_{1, N}(-m,n)$ in the holomorphic part of $\mathcal{P}_{1,N}(-m,z)$ is given by
 $$b_{1, N}(-m,n)= -\pi \sqrt 2 (n/m)^{1/4} (1+i) \sum_{c>0 \atop 4N|c} \delta_{o} (c/4)
   \frac{K_1(-m,n,c)}{c} I_{1/2}\left( \frac{4\pi\sqrt{nm}}{c} \right). $$
 Applying \cite[Proposition 3.1]{BO}, we may write the right-hand side as
 $$-\pi \sqrt 2 (n/m)^{1/4} (1+i) \sum_{c>0 \atop 4N|c} \delta_{o} (c/4)(-i)
   \frac{K_0(-m,n,c)}{c} I_{1/2}\left( \frac{4\pi\sqrt{nm}}{c} \right),$$
 which is the Fourier coefficients $-b_{0,N}(-n,m)$ in the holomorphic part of $\mathcal{P}_{0,N}(-m,z)$ by \cite[Theorem 2.1]{MP}.
Thus
 \begin{equation} \label{D1} b_{1,N}(-m,n)=-b_{0,N}(-n,m).\end{equation}
By the definition of $B_{\lambda, N}(-m,n)$, we have
 \begin{equation} \label{D2}
 B_{1,N}(-m,n)=-c_1 \delta_{\square,m} H_N^* (n)+ b_{1,N}(-m,n)
 \end{equation} and
 \begin{equation} \label{D3}
 B_{0,N}(-m,n)=c_1 \delta_{\square,n} H_N^* (m)+ b_{0,N}(-m,n).
 \end{equation}
 Now the theorem follows from (\ref{D1}), (\ref{D2}) and (\ref{D3}).
\end{proof}

\section{Congruences of traces}

We start the proof of Theorem \ref{last} by showing that the Niebur-Poincar\'e series $\mathcal F_m (z,1)$ is generated by the action of Hecke operator on $\mathcal F_1 (z,1)$.

\begin{lem} \label{Hecke} For every positive integer $m$ relatively prime to $N$, we have
\begin{equation}\label{mhe} \mathcal F_1^*(z,1)|_{T_m}= \mathcal{F}_m^*(z,1),
\end{equation}
where $|_{T_m}$ denotes the usual action of the $m$th Hecke operator $T_m$ on the space of weakly holomorphic modular functions.
And for a prime $p$ not dividing $N$,
\begin{equation}\label{hecke} \mathcal F_m^*(z,1)|_{T_p}= \mathcal{F}_{pm}^*(z,1)+p\mathcal{F}_{m/p}^*(z,1),
\end{equation} where $\mathcal{F}_{m/p}^*(z,1)$ is defined to be zero unless $m/p \in \mathbb Z$.
\end{lem}
\begin{proof} As we can easily deduce (\ref{mhe}) from (\ref{hecke}), we only prove (\ref{hecke}).  By \cite[Lemma 6]{AL} and \cite[Lemmas 2.5 and 2.6]{Koike}, $\mathcal F_m^*(z,1)|_{T_p}$ is on $\Gamma_0^*(N)$ and has a pole only at
 $\infty$. Since $\mathcal F_m^*(z,1)$ is harmonic, so is  $\mathcal F_m^*(z,1)|_{T_p}$. Now if we compare the principal parts of the functions in both sides of (\ref{hecke}) and apply \cite[Theorem 6]{Niebur}, we obtain the result immediately.
\end{proof}

By Theorem \ref{Kloo}, for $m$ and $p$ satisfying $(m,N)=1$ and $(p,N)=1$, we have
\begin{eqnarray} \label{LHS}
\textbf{t}_{pm}(D)&=&-\left(\sum_{\nu\mid m}\nu B(-\nu^2,D)+\sum_{\nu\mid m}p\nu B(-p^2\nu^2,D) \right)\cr
&=& \textbf{t}_{m}(D)-p\sum_{\nu\mid m}\nu B(-p^2\nu^2,D).
\end{eqnarray}
On the other hand, by the definition of $\textbf{t}_{m}(D)$ and (\ref{hecke}) (cf. \cite[Section 6]{Zagier}),
\begin{equation} \label{RHS}
\textbf{t}_{pm}(D)=\left(\textbf{t}_{m}(p^2 D)+ \left( \frac{-D}{p}\right) \textbf{t}_{m}(D)+p \textbf{t}_{m}(D/p^2)\right)-p \textbf{t}_{m/p}(D).
\end{equation}
Since $\textbf{t}_{f}(D)=\sum_ma_m\textbf{t}_{m}$, it follows from (\ref{LHS}) and (\ref{RHS}) that
\begin{equation} \label{star}
\textbf{t}_{f}(p^2 D)+ \left( \frac{-D}{p}\right) \textbf{t}_{f}(D)+p \textbf{t}_{f}(D/p^2)=\textbf{t}_{f}(D)-p C(D)
\end{equation}
where \begin{equation*}C(D)= \sum_ma_m
\left(\sum_{\nu\mid m}\nu B(-p^2\nu^2,D)+\sum_{\nu\mid m/p}\nu B(-\nu^2,D)\right).\end{equation*}
Therefore, for a prime $p$ for which $\left( \frac{-D}{p}\right)=1 $, we find that
\begin{equation} \label{short}
\textbf{t}_{f}(p^2 D)=- p C(D)
\end{equation}
which proves (\ref{cons1}).

\section{The limiting distribution of the modular traces}

Using his result on uniform distribution of CM points on the modular curve $X_0(1)$ \cite{Duke2}, Duke \cite[Theorem 1]{Duke} proved a conjecture on the limiting distribution of traces of singular moduli $j_1(z)$ made by Bruinier, Jenkins, and Ono \cite[Theorem 1.1]{BJO}. The limiting distribution of traces of singular moduli, so-called \lq\lq 24 Theorem\rq\rq, is equivalent to
 \begin{equation}\label{24}
 \mathop{\sum_{c>\sqrt{D/3}}}_{c\equiv 0
(\textrm{mod}\ 4)}
 S_D(c)\sinh(4\pi\sqrt{D}/c)=o(H(D)),\end{equation}
where the Sali\'{e} sum $S_D(c)$ is defined for any positive integer $c$ by $\displaystyle{S_D(c):=\sum_{x^2\equiv -D (\textrm{mod}\, c)} e(2x/c)}$.
This has been recently generalized to the limiting distribution of the twisted trace of Maass-poincar\'{e} series  $\mathcal F_m (z,s)$ for $\G(1)$ by Folsom and Masri \cite[Theorem 1.1]{FM}.  Although we could evaluate exact values of $\sum_c
 S_D(m,c)\sinh\left(\frac{4\pi m\sqrt{D}}{c}\right)$ for many cases as in Examples 1 and 2 in Section \ref{results}, we were not able to establish the analogue of the equation (\ref{24}) due to the difficulty of determining a proper bound in the fundamental domain of $\G_0(N)$ which contains enough of Heegner points in $X_0(N)$. From our computational experiment in level 4 case, though, we find that
\begin{equation}\label{j4}
 \mathop{\sum_{c>16\sqrt{D/3}}}_{c\equiv 0
(\textrm{mod}\ 4)}
 S_D(c)\sinh(4\pi\sqrt{D}/c)=o(H(D)).\end{equation}


\begin{thebibliography}{99}
\bibitem{AO}
S.~Ahlgren and K.~Ono, \emph{Arithmetic of Singular Moduli and Class Equations}, Compositio Math. \textbf{141} (2005), 293-312.
%
\bibitem{AL}
A.~O.~L.~Atkin and J.~Lehner,
\emph{Hecke operators on $\Gamma_0(m)$},
Math. Ann. 185 (1970), 134-160.
%

\bibitem{Borch}
R.~E.~Borcherds, \emph{Automorphic forms on $O_{s+2,2}(\mathbb{R})$ and infinite products}, Invent. Math. \textbf{120} (1) (1995), 161--213.
%

\bibitem{Bru}
J.~H.~Bruinier, \emph{Borcherds products on $O(2,\ell)$ and Chern classes of Heegner divisors}, Springer Lect. Notes, 1780, Springer-Verlag, Berlin, 2002.
%

\bibitem{BF}
J.~H.~Bruinier and J.~Funke, \emph{Traces of CM-values of modular
functions}, J.~Reine Angew.~Math.~\textbf{594} (2006), 1--33.
%

\bibitem{BJO}
J.~H.~Bruinier, P.~Jenkins, and K.~Ono, \emph{Hilbert class
polynomials and traces of singular moduli},  Math.~Ann.~
\textbf{334} (2006), 373--393.
%
\bibitem{BO}
K.~Bringmann and K.~Ono,
\emph{Arithmetic properties of coefficients of half-integral weight Maass-Poincar\'e series}, Math.~Ann.~
\textbf{337} (2007), 591--612.


\bibitem{CJKK1}
D.~Choi, D.~Jeon, S.-Y.~Kang, and C.~H.~Kim,
\emph{Traces of singular moduli of arbitrary level modular functions}, IMRN \textbf{2007} (2007), article ID rnm110, 17 pages.
%
\bibitem{CJKK2}
D.~Choi, D.~Jeon, S.-Y.~Kang, and C.~H.~Kim,
\emph{Exact formulas for traces of singular moduli of higher level modular functions}, J. Number T. \textbf{128} (3) (2008), 700--707.
%
\bibitem{Duke}
W.~Duke, \emph{Modular functions and the uniform distribution of CM
 points},  Math.~Ann.~ \textbf{334} (2006), 241--252.
%
\bibitem{Duke2}
W.~Duke, \emph{Hyperbolic distribution problems and half-integral weight Maass forms},  Invent.~Math.~ \textbf{92} (1988), 73--90.
%
\bibitem{FM}
A.~Folsom and R.~Masri, \emph{Limiting distribution of traces of Maass-Poincar\'{e} series}, preprint (http://www.math.wisc.edu/$\sim$masri/Maass-Poincare-IMRN.pdf).
%
\bibitem{Hej}
D.~A.~Hejhal, \emph{The Selberg trace formula for $PSL_2(\mathbb{R})$}, Springer Lect. Notes, 1001, Springer-Verlag, Berlin, 1983.
%
\bibitem{Koike}
M. Koike,
\emph{On Replication formula and Hecke operators},
preprint.
%
\bibitem{Kohnen}
W.~Kohnen, \emph{Fourier coefficients of modular forms of half-integral weight},
Math. Ann. \textbf{271} (1985), 237--268.
%
\bibitem{GKZ}
{B.~Gross, W.~Kohnen, and D.~Zagier}, \emph{Heegner points and derivatives of $L$-seires, II}, Math. Ann. \textbf{278} (1987), 497--562.
%
\bibitem{MP}
A.~Miller and A.~Pixton,
\emph{Arithmetic traces of non-holomorphic modular invariants},
preprint (http://www.math.wisc.edu/$\sim$ono/reu06zagier.pdf).
%
\bibitem{Niebur}
D.~Niebur, \emph{A class of nonanalytic automorphic functions},
Nagoya Math.~ J.~ \textbf{52} (1973), 133--145.
%
\bibitem{ono}
K.~Ono, \emph{Unearthing The Vision of a Master}, Harvard-MIT Current Developments in Mathematics 2008, International Press, (preliminary version: http://www.math.wisc.edu/$\sim$ono/reprints/114.pdf).
%
\bibitem{Os}
R.~Osburn, \emph{Congruences for Traces of Singular moduli}, Ramanujan J. \textbf{14} (2007) 411--419.
%
\bibitem{Zagier}
D.~Zagier, \emph{Traces of singular moduli}, Motives,
polylogarithms and Hodge theory, Part I (Irvine, CA, 1998),
211--244, Int. Press Lect. Ser., 3, I, Int. Press, Somerville, MA,
2002.

\end{thebibliography}
\end{document}